\documentclass[12pt,reqno]{amsart}
\usepackage{amssymb,a4wide,amsthm,amsmath,amsfonts}
\usepackage{mathrsfs}
\def\mathcal{\mathscr}
\usepackage[colorlinks,linkcolor={blue},
citecolor={blue},urlcolor={red},]{hyperref}
\usepackage{hyperref}
\let\emptyset \undefined
\let\ge       \undefined
\let\le       \undefined
\newsymbol\le          1336  
\newsymbol\ge          133E  
\newsymbol\emptyset    203F
\newsymbol\notle       230A
\newsymbol\notge       230B
\theoremstyle{plain}
\newtheorem{theorem}{Theorem}[section]
\theoremstyle{remark}
\newtheorem{remark}[theorem]{Remark}
\newtheorem{example}[theorem]{Example}
\theoremstyle{plain}
\newtheorem{corollary}[theorem]{Corollary}
\newtheorem{lemma}[theorem]{Lemma}
\newtheorem{proposition}[theorem]{Proposition}
\newtheorem{definition}[theorem]{Definition}
\newtheorem{hypothesis}[theorem]{Hypothesis}
\numberwithin{equation}{section}
\begin{document}
\title[]
{THE ORNSTEIN UHLENBECK BRIDGE AND APPLICATIONS TO MARKOV SEMIGROUPS}
\author{B. Goldys}
\address{School of Mathematics, The University of New South
Wales, Sydney 2052, Australia}
\email{B.Goldys@unsw.edu.au}
\author{B. Maslowski}
\address{Institute of Mathematics, Academy of Sciences of Czech Republic\\
 \v Zitn\' a 25, 11567 Praha 1, Czech Republic}
\email{maslow@math.cas.cz}
\thanks{This work was partially supported by the UNSW Faculty Research Grant 
 and GA\v CR grant 201/04/0750} \keywords{Ornstein-Uhlenbeck process, pinned process, measurable linear mapping, stochastic semilinear equation,
transition density}
\subjclass{60H15, 60H10, 60J60, 60J35, 35R60}
\begin{abstract}
For an arbitrary Hilbert space-valued Ornstein-Uhlenbeck process we construct the Ornstein-Uhlenbeck Bridge connecting
a starting point $x$ and an endpoint $y$ that belongs to a certain linear subspace of full measure. We derive also a
stochastic evolution equation satisfied by the OU Bridge and study its basic properties. The OU Bridge is then used to
investigate the Markov transition semigroup associated to a nonlinear stochastic evolution equation with additive noise. We provide
an explicit formula for the transition density and study its regularity. Given the Strong Feller property and the
existence of an invariant measure we show that the transition semigroup maps $L^p$ functions into continuous functions.
We also show that transition operators are $q$-summing for some $q>p>1$, in particular of Hilbert-Schmidt type.

\end{abstract}
\maketitle \tableofcontents
 \section{Introduction}
Let $\left(Z_t^x\right)$ be an Ornstein-Uhlenbeck process on a
separable Hilbert space $H$. By this we mean that $\left(Z_t^x\right
)$ is
a solution to a linear stochastic evolution equation
\begin{equation}\left\{\begin{array}{l}
dZ_t^x=AZ_t^xdt+\sqrt {Q}dW_t,\\
Z_0^x=x\in H.\end{array}
\right.\label{01}\end{equation}
In the above equation $\left(W_t\right)$ is a standard cylindrical
Wiener process defined on a certain stochastic basis
$\left(\Omega ,\mathcal{F},\left(\mathcal{F}_t\right),\mathbb{P}\right
)$ and $Q=Q^{*}\ge 0$ is a
bounded operator on $H$. We assume that the operator $\left(A,\mbox{\rm dom}
(A)\right)$ is a
generator of a $C_0$-semigroup $\left(S_t\right)$ on $H$. Under the assumptions given below the solution to (\ref{01}) is defined by the formula
\begin{equation}Z_t^x=S_tx+\int_0^tS_{t-s}\sqrt QdW_s.\label{02}\end{equation}
The aim of this paper is to study the basic properties of the Ornstein-Uhlenbeck Bridge (sometimes called a Pinned Ornstein-Uhlenbeck process)
$\left(\hat Z_t^{x,y}\right)$ associated to the Ornstein-Uhlenbeck process $\left(Z_t^x\right)$ and its applications. Let us recall informally, that
this process is defined via the formula
\[\mathbb P\left(\left.Z_t^x\in B\right|Z_T^x=y\right)=\mathbb P\left
(\hat Z_t^{x,y}\in B\right),\quad t<T,\] where $x,y\in H$ and $B\subset H$ is a Borel set. Intuitively, it is an
Ornstein-Uhlenbeck process "conditioned to go from $x$ at time $t=0$ to $y$ at time $t=T$" (a rigorous definition is
given in Section 2, cf. Def. \ref{OUB}). The importance of various types of bridge processes in the theory of finite
dimensional diffusions is well recognised, see for example \cite{yor}. In infinite dimensional framework this concept
was developed in \cite{simao1} in order to study regularity of transition semigroup of certain linear and nonlinear
diffusions on Hilbert space. In \cite{masi1} and \cite{masi2} an Ornstein-Uhlenbeck Bridge is introduced in order to
obtain lower estimates on the transition kernel of some semilinear stochastic evolution equations. Those estimates
provide a powerful tool to study exponential ergodicity and $V$-uniform ergodicity for such equations. In particular,
they allowed us to obtain in \cite{den} explicit estimates of the rate of exponential convergence to the invariant
measure.

In the present paper the OU Bridge is studied under much more general conditions and in more detail. We provide also
further applications of the OU Bridge to the analysis of transition densities and the regularity of associated Markov
semigroups. Regularity of Strongly Feller transition semigroups was studied in \cite{furman} (see also references
therein). We use methods completely different from \cite{furman} and obtain stronger results but for bounded drifts
only while the aforementioned paper allows linearly growing drifts. Closely related results for semigroups that are not
strongly Feller may be found in \cite{ania}. For the regularity of strongly Feller semigroups associated to the OU
process we refer to \cite{reg}.

Let us describe the contents of this paper. In Section 2 we provide, for the reader's convenience, some relevant facts
about linear measurable mappings and conditional distributions of Hilbert space valued Gaussian random vectors. Then we
give a definition of the OU Bridge and some basic results on OU processes and OU Bridges. Some of the technical results
from \cite{den} that are needed in the sequel are stated without proof and others (Lemma \ref{vt1a}, Proposition
\ref{ou0} and Lemma \ref{lh1}) are reproved under more general conditions. In Section 3, a stochastic equation for the
OU Bridge is derived. A new Brownian Motion adapted to the filtration of the Ornstein Uhlenbeck Bridge is obtained and
then it is shown that the Bridge process is a unique mild (and weak) solution of a linear nonhomogeneous stochastic
evolution equation with singular coefficients. Section 4 is devoted to applications of the previous results to
semilinear stochastic equations; continuity of Markov transition densities (with respect to the Gaussian invariant
measure $\nu $ that is an invariant measure with respect to the OU process) is proved (Theorem \ref{dens} and Remark
\ref{iny}), the Markov semigroup is shown to map the space $L^p(H,\nu )$, $p>1$, into the space of continuous functions
on $H$ (Theorem \ref{lp}) and is also shown to be Hilbert-Schmidt on $L^2(H, \nu )$ and $q$-summing (in particular,
compact) as a mapping $L^p(H, \nu ) \to L^q(H, \nu )$ even if $q>p$ provided the gap between $q$ and $p$ is not too
large (Theorem \ref{HSS}). At the end of the section the results are illustrated in the case of one-dimensional
semilinear stochastic parabolic equation (Example \ref{example}) in which case all conditions imposed in the paper are
verified or specified.

ACKNOWLEDGEMENT. The authors are grateful to Jan Seidler for his valuable comments and suggestions.

\section{Preliminaries on OU Processes and Bridges}
In this section we collect, for the reader convenience,
some properties of infinite-dimensional OU processes and
Gaussian random variables which will be useful in the
paper. We also define the OU Bridge and recall some known results that will be useful in the sequel.
\subsection{Measurable Linear Mappings}
Let $H$ be a real separable Hilbert space and let $\mu =N(0,C)$ be a centered Gaussian measure on $H$ with the covariance operator $C$ such that
$\overline {\mbox{\rm im}(C)}=H$. The space $H_C=\mathrm{im}\left(C^{1/2}\right)$ endowed with the norm $|x|_ C=\left|C^{-1/2}x\right|$ can be
identified as the Reproducing Kernel Hilbert Space of the measure $\mu$. In the sequel we will denote by $\left\{e_n:n\ge 1\right \}$ the eigenbasis
of $C$ and by $\left\{c_n:n\ge 1\right\}$ the corresponding set of eigenvalues:
\[Ce_n=c_ne_n,\quad n\ge 1.\]
For any $h\in H$ we define
\[\phi_n(x)=\sum_{k=1}^n\frac 1{\sqrt {c_k}}\left\langle h,e_k\right
\rangle\left\langle x,e_k\right\rangle ,\quad x\in H.\]
The following two lemmas are well known (see e.g. \cite{den}):
\begin{lemma}\label{fiha}
The sequence $\left(\phi_n\right)$ converges in $L^2(H,\mu )$ to a limit $
\phi$
and
\[\int_H\left|\phi (x)\right|^2\mu (dx)=|h|^2.\]
Moreover, there exists a measurable linear space
$\mathcal{M}_h\subset H$, such that $\mu\left(\mathcal{M}_h\right)=
1$, $\phi$ is
linear on $\mathcal M_h$ and
\begin{equation}\phi (x)=\lim_{n\to\infty}\phi_n(x),\quad x\in\mathcal
M_h.\label{limit}\end{equation}
We will use the notation $\phi (x)=\left\langle h,C^{-1/2}x\right
\rangle$.
\end{lemma}
Let $H_1$ be another real separable Hilbert space and let
$T:H\to H_1$ be a bounded operator. The Hilbert-Schmidt norm
of $T$ will be denoted by $\left\|T\right\|_{HS}$. Let
\[\tilde {T}_nx=\sum_{k=1}^n\frac 1{\sqrt {c_k}}\left\langle x,e_
k\right\rangle Te_k,\quad x\in H.\]
\begin{lemma}\label{tmeas}
Let $T:H\to H_1$ be a Hilbert-Schmidt operator. Then the
sequence $\left(\tilde T_n\right)$ converges in $L^2\left(H,\mu ;
H_1\right)$ to a limit $\tilde {T}$
and
\[\int_H\left|\tilde T(x)\right|^2_{H_1}\mu (dx)=\left\|T\right\|_{
HS}^2.\]
Moreover, there exists a measurable linear space
$\mathcal M_T\subset H$, such that $\mu\left(\mathcal M_T\right)=
1$, $\tilde {T}$ is
linear on $\mathcal M_T$ and
\begin{equation}\tilde {T}(x)=\lim_{n\to\infty}\tilde {T}_nx,\quad
x\in\mathcal M_T.\label{limit1}\end{equation}
We will use the notation $TC^{-1/2}x := \tilde T(x)$.
\end{lemma}
The above procedure is specified in the following Lemma (the proof of which may be found in \cite{den}) to operator-valued functions:
\begin{lemma}\label{integral}
Let $K(t,s):H\to H$ be an operator-valued, strongly measurable function, such that for each $a\in (0,T)$
\begin{equation}\int_0^a\int_0^a\left\|K(t,s)\right\|^2_{HS}dsd
t<\infty .\label{double}\end{equation} Then the following holds. \par\noindent
(a) There exists a Borel set $B\subset [0,T]^2$ of full Lebesgue
measure such that the measurable linear mapping $K(t,s)C^{-1/2}$ is well defined for all $(s,t)\in B$.
\par\noindent
(b) There exists a measurable mapping $f:[0,T)^2\times H\to H$ and a measurable linear space $\mathcal M\subset H$ of full measure such that
$f(t,s,y)=K(t,s)C^{-1/2}y$ for $y\in \mathcal M$ and for each $a<T$
\[\int_0^a\left|f(t,s,y)\right|ds<\infty \]
for almost all $t\in [0,T]$. We will use the notation $K(t,s)C^{-1/2}y := f(t,s,y)$.
\end{lemma}
\vskip1cm
\subsection{Conditional Distributions}
\par\noindent
Let $H_1$ and $H_2$ be two real separable Hilbert spaces and
let $(X,Y)\in H_1\times H_2$ be a Gaussian vector with mean values
\[m_X=\mathbb EX,\quad\mbox{\rm and}\quad m_Y=\mathbb EY.\]
 The covariance
operator of $X$ is determined by the equation
\begin{equation}\mathbb E\left\langle X-m_X,h\right\rangle\left\langle
X-m_X,k\right\rangle =\left\langle C_Xh,k\right\rangle ,\quad h,k
\in H_1,\label{con1}\end{equation}
and a similar condition determines the covariance $C_Y$ of
$Y$. The covariance operator $C_{XY}:H_1\to H_2$ is defined by the
condition
\[\left\langle C_{XY}h,k\right\rangle =\mathbb E\left\langle X-m_
X,h\right\rangle\left\langle Y-m_Y,k\right\rangle ,\quad h\in H_1 ,k\in H_2,\] and then $C_{XY}^{*}=C_{YX}$.\\
 For a
linear closable operator $G$ on $H$ the closure of $G$ will be denoted by $\overline {G}$. The next theorem is well
known, see for example \cite{mandelbaum}
\begin{theorem}\label{tcond}
Assume that $C_X$ is injective.
Then the following holds.
\par\noindent
(a) We have
\begin{equation}\mbox{\rm im}\left(C_{YX}\right)\subset\mbox{\rm im}\left
(C_X^{1/2}\right),\label{con2}\end{equation}
the operator
$T=C_X^{-1/2}C_{YX}$ is of
Hilbert-Schmidt type on $H$ and $T^{*}=\overline {C_{XY}C_X^{-1/2}}$.
\par\noindent
(b) We have
\[\mathbb E\left(Y|X\right)=m_Y+T^{*}C_X^{-1/2}\left(X-m_X\right)
,\quad\mathbb P^X-a.s.\]
\par\noindent
(c) The conditional distribution of $Y$ given $X$ is Gaussian
$N\left(\mathbb E\left(Y|X\right),C_{Y|X}\right)$, where
\[C_{Y|X}=C_Y-T^{*}T.\]
Moreover, the random variables $T^{*}C_X^{-1/2}X$ and
$\left(Y-T^{*}C_X^{-1/2}X\right)$ are independent.
\end{theorem}
\vskip1cm
\subsection{Some Properties of the Ornstein-Uhlenbeck Process}

\par\noindent
The following hypothesis is a standing assumption for the rest of the paper.
\begin{hypothesis}\label{hou}
For every $t>0$
\begin{equation}\int_0^t\left\|S_sQ^{1/2}\right\|^2_{HS}ds<\infty
,\label{hs1}\end{equation} and
\begin{equation}
\overline {\mbox{\rm im}\left(Q_t\right)}=H,\label{hs2}
\end{equation}
where, in view of (\ref{hs1})
\begin{equation}\label{ouqt}
Q_t=\int_0^tS_sQS_s^*ds.
\end{equation}
is a well defined trace class operator.
\end{hypothesis}
\noindent It is well known that if Hypothesis \ref{hou} holds then the process (\ref{02}) is a well defined $H$-valued,
Gaussian and Markov process, see \cite{dz1}.
\par
Let $\mu$ denote the probability law of the process $\left\{Z_t^0:t\in [0
,1]\right\}$ that is
concentrated on $L^2(0,T;H)$ and let
$\mathcal L:L^2(0,T;H)\to C(0,T;H)$
be defined by the formula
\begin{equation}\label{rkhs}
\mathcal Lu(t)=\int_0^tS_{t-s}Q^{1/2}u(s)ds.
\end{equation}
Note that, cf. \cite{dz1}, $\mathrm{im}(\mathcal L)=RKHS(\mu )$  (the Reproducing Kernel Hilbert Space of the measure
$\mu$). We will use the notation $\mu_t^x$ for the Gaussian measure $N\left (S_tx,Q_t\right)$ and $\mu_t$ for
$\mu_t^0$. By the properties of Gaussian distribution $\mu_t^x$ is the probability distribution of a random variable
$Z_t^x$ and we set $Z_t = Z^0_t$. In the rest of this subsection we give several statements on properties of the family
of covariance operators $\left\{ Q_t:\,  t\le T\right\} $ that will be useful later.
\par\noindent
The definition of $Q_t$ given in (\ref{ouqt}) yields immediately a simple identity that will be frequently used:
\begin{equation}\label{qtd}
Q_T=Q_t+S_tQ_{T-t}S_t^*,\quad t\le T.
\end{equation}
\begin{lemma}\label{lqt}
We have
\[\mathrm{im}\left(Q_t^{1/2}\right)\subset\mathrm{im}\left(Q_T^{1/2}\right),\quad t\le T,\]
hence the operator $U_t=Q_T^{-1/2}Q_t^{1/2}$ is bounded on $H$ for every $t\le T$ and  $\left\|U_t\right\|\le 1$.
Moreover, $U_t^*=\overline{Q_t^{1/2}Q_T^{-1/2}}$, the closure of the operator $Q_t^{1/2}Q_T^{-1/2}$ defined on the
domain $\mathrm{im}\left(Q_T^{1/2}\right)$.
\end{lemma}
\begin{proof}
From the definition of the covariance operators $Q_t$ it follows that $|Q_tx|^2 \le |Q_Tx|^2$ for each $x\in H$ and
$0\le t\le T$ and the conclusion easily follows.
\end{proof}
\begin{lemma}\label{vt1}
(a) The operator
$V_t=Q_T^{-1/2}S_{T-t}Q_t^{1/2}$ is well defined and bounded on $H$ and
\begin{equation}\left\|V_t\right\|\le 1,\quad t\in (0,T).\label{nonexp}\end{equation}
Moreover,
\begin{equation}\lim_{t\to T}V_t^{*}x=\lim_{t\to T}V_tx=x,\quad x
\in H.\label{c0}\end{equation}
(b) For any $t\in [0,T]$
\begin{equation}Q_{T-t}=Q_T^{1/2}\left(I-V_tV_t^{*}\right)Q_T^{1/
2}.\label{contr1}\end{equation}
\end{lemma}
\begin{proof}
The inequality (\ref{nonexp}) has been proved in \cite{neerven}, the convergence (\ref{c0}) in \cite{den}.
Part (b) follows immediately from (\ref{qtd}).
\end{proof}
Under a slightly stronger condition we show that the inequality (\ref{nonexp}) is sharp, more precisely, we have
\begin{lemma}\label{vt1a}
The following conditions are equivalent:
\par\noindent
(a) For any $t\in (0,T]$
\begin{equation}\label{hh2}
\mbox{\rm im}\left(Q_t^{1/2}\right)=\mbox{\rm im}\left(Q_T^{1/2}\right
).
\end{equation}
(b) $\mathrm{im}\left(U_t\right)$ is dense in $H$ for each $t\in (0,T)$.
\par\noindent
(c) We have
\begin{equation}\label{contr}
\left\|V_t\right\|  <1,\quad t \in (0,T).
\end{equation}
\end{lemma}
\begin{proof}
Obviously (a) implies (b).
\par\noindent
 To prove that (b) implies (c) note first that  (\ref{qtd}) yields
\[\left|Q_{T-t}^{1/2}x\right|^2=\left|Q_T^{1/2}x\right|^2-\left|V_t^*Q_T^{1/2}x\right|^2,\]
hence putting $y=Q_T^{1/2}x$ we obtain
\[\left|Q_{T-t}^{1/2}Q_T^{-1/2}y\right|^2=\left|y\right|^2-\left|V_t^*y\right|^2.\]
Assume that $\left\|V_t^*\right\|=1$ for a certain $t\in(0,T)$. Since $\mathrm{im}\left(Q_T^{1/2}\right)$ is dense in $H$, there exists a sequence
$y_n\in\mathrm{im}\left(Q_T^{1/2}\right)$, such that $\left|y_n\right|=1$ and $\left|V_t^*y_n\right|\to 1$. Therefore,
\begin{equation}\label{lim0}
\lim_{n\to\infty}\left|Q_{T-t}^{1/2}Q_T^{-1/2}y_n\right|^2=\lim_{n\to\infty}\left(1-\left|V_t^*y_n\right|^2\right)=0.
\end{equation}
 Let $y_{n_k}$ be a subsequence
converging weakly to $y\in H$. Since
\[\mathrm{im}\left(Q_{T-t}^{1/2}\right)\subset\mathrm{im}\left(Q_T^{1/2}\right),\quad t\le T,\]
and
\[\left(Q_T^{-1/2}Q_{T-t}^{1/2}\right)^*=\overline{Q_{T-t}^{1/2}Q_T^{-1/2}},\]
we find that
\[Q_{T-t}^{1/2}Q_T^{-1/2}y_{n_k}\to \overline{Q_{T-t}^{1/2}Q_T^{-1/2}}y,\quad\mathrm{weakly},\]
and by (\ref{lim0}) we obtain $\overline{Q_{T-t}^{1/2}Q_T^{-1/2}}y=0$ and since $\left|V_t^*y\right|=1$ we obtain $y\neq 0$. It follows that the
range of the operator $Q_T^{-1/2}Q_{T-t}^{1/2}$ is not dense in $H$, which shows that (b) implies (c).
\par\noindent
Finally, assume that (c) holds. Then (\ref{contr1}) and Proposition B1 in \cite{dz1} yield
\[\mathrm{im}\left(Q_{T-t}^{1/2}\right)=\mathrm{im}\left(Q_T^{1/2}\left(I-V_tV_t^*\right)^{1/2}\right).\]
Since $\left\|V_t\right\|<1$, the operator $I-V_tV_t^*:H\to H$ is an isomorphism, hence
\[\mathrm{im}\left(Q_{T-t}^{1/2}\right)=\mathrm{im}\left(Q_T^{1/2}\right),\quad t<T,\]
and (a) follows.
\end{proof}
\begin{remark}\label{remh2}
Necessary and sufficient conditions for  (\ref{hh2}) to hold are not known but it was proved to be satisfied in the following cases.
 \par\noindent
 (a) If
 \[\mathrm{im}\left(S_t\right)\subset\mathrm{im}\left(Q_t^{1/2}\right), \quad t>0,\]
  then (\ref{hh2}) holds. It is known that the above condition is equivalent to the strong Feller property of the OU transition semigroup $R_t\phi(x)=\mathbb E\phi\left(Z_t^x\right)$, see \cite{dz1} for details.
 \par\noindent
 (b) Assume that the process $\left(Z_t^x\right)$ admits a nondegenerate invariant measure $\nu$ and $\mathrm{im}(Q)$ is dense in $H$. Let $H_Q=\mathrm{im}\left(Q^{1/2}\right)$ be endowed with the norm $|x|_Q=\left|Q^{-1/2}x\right|$. Assume that $H_Q$ is invariant for the semigroup $\left(S_t\right)$ and its restriction to $H_Q$ is a $C_0$-semigroup in $H_Q$. Then (\ref{hh2}) holds, see \cite{acta}. These assumptions are satisfied for any process $\left(Z_t^x\right)$ with the transition semigroup analytic in $L^2(H,\nu )$, in particular they are satisfied for any reversible OU process.
\end{remark}
We define the operator $B: Q_T^{1/2}(H) \to L^2(0,T;H)$,
\[Bx(t)=Q^{1/2}S_{T-t}^{*}Q_T^{-1/2}x,\quad t\in [0,T],\quad x\in Q_T^{1/2}(H).\]
\par\noindent
 The following simple Lemma has been proved in \cite{den}:
\begin{lemma}\label{ania}
(a) The operator $B$ with the domain $\mathrm{dom}(B)=Q_T^{1/2}(H)$ extends to a bounded operator (still denoted by
$B$) $B:H\to L^2( 0,T;H)$. Moreover,
\[|Bx|_{L^2(0,T;H)}=|x|_H,\quad x\in H.\]
\par\noindent
(b) Seting
\begin{equation}
H\ni x\to\mathcal Kx(t)=K_tx\in L^2(0,T;H),\label{dk}
\end{equation}
where
\begin{equation}K_t=Q_t^{1/2}V_t^{*},\label{kt}\end{equation}
we have $\mathcal K=\mathcal LB$.
In particular the operator $\mathcal K:H\to C(0,T;H)$ is bounded.
\end{lemma}
\subsection{Fundamentals on OU Bridge}
In the present subsection we give the definition and some basic
properties of the OU Bridge.
\par\noindent
Since $V_t^{*}=\overline {Q_t^{1/2}S_{T-t}^{*}Q_T^{-1/2}}$ is bounded, the operator $K_t$ is of Hilbert-Schmidt type on
$H$ for each $t\in [0,T)$. Also, $\mathcal K : H\to L^2(0,T;H)$ is Hilbert-Schmidt.\\
 Note that if $K_t$ is defined by
(\ref{kt}) then, in view of Lemma \ref{tmeas}, the measurable function $K_tQ_T^{-1/2}$ is well defined for each $t\in
[0,T]$. We will start from the definition  of the process $(\hat Z_t )$,
\[\hat {Z}_t=Z_t-K_tQ_T^{-1/2}Z_T,\quad t\in [0,1),\quad\mbox{\rm and}
\quad\hat {Z}_1=0.\]
\begin{proposition}\label{ou0}
(a) An $H$-valued Gaussian process $\left(\hat Z_t\right)$ is independent of $
Z_T$.
\par\noindent
(b) The covariance operator $\hat {Q}_t$ of $\hat {Z}_t$ is given by
\begin{equation}\hat {Q}_t=Q_t^{1/2}\left(I-V_t^{*}V_t\right)Q_t^{
1/2}.\label{covp}\end{equation}
(c) The process $\left(\hat Z_t\right)$ is mean-square continuous on $[0,T]$.
\par\noindent
(d) If, moreover, one of the equivalent conditions (a)-(c) of Lemma \ref{vt1a} holds then
\begin{equation}\mbox{\rm im}\left(\hat Q_t^{1/2}\right)=\mbox{\rm im}\left
(Q_t^{1/2}\right),\quad t\in (0,T).\label{im12}\end{equation}
\end{proposition}
\begin{proof}
Theorem \ref{tcond} yields immediately (a) since
$\hat {Z}_t=Z_t - \mathbb E\left(\left.Z_t\right|Z_T\right)$. Invoking (c) of Theorem \ref{tcond}\
with $C_X=Q_T$, $C_Y=Q_t$ and $T^{*}=K_t$ and (\ref{kt}) we obtain
\[\hat {Q}_t=Q_t-K_tK_t^{*}=Q_t^{1/2}\left(I-V_t^{*}V_t\right)Q_t^{
1/2},\quad t<T.\] Using (\ref{nonexp}) we find easily that
\begin{equation}\lim_{t\to 0}\mbox{\rm tr}\left(\hat Q_t\right)=0
.\label{c00}\end{equation}
To prove that
\begin{equation}\lim_{t\to T}\mbox{\rm tr}\left(\hat Q_t\right)=0
,\label{c1}\end{equation}
we note first that
\[\mbox{\rm tr}\left(\hat Q_t\right)=\mbox{\rm tr}\left(\left(I-V_
t^{*}V_t\right)\left(Q_t-Q_T\right)\right)+\mbox{\rm tr}\left(\left
(I-V_t^{*}V_t\right)Q_T\right).\]
Next, it is easy to see that
\begin{equation}0\le\lim_{t\to T}\mbox{\rm tr}\left(\left(I-V_t^{
*}V_t\right)\left(Q_T-Q_t\right)\right)\le\lim_{t\to T}\mbox{\rm tr}\left
(Q_T-Q_t\right)=0.\label{z0}\end{equation}
Finally,
\[\mbox{\rm tr}\left(\left(I-V_t^{*}V_t\right)Q_T\right)=\mbox{\rm tr}\left
(Q_T\right)-\mbox{\rm tr}\left(V_tQ_TV_t^{*}\right)\]
\[=\mathrm{tr}\left(Q_T\right)-\sum_{k=1}^{\infty}\left|Q_T^{1/2}
V_t^{*}e_k\right|^2,\]
where $\left\{e_k:k\ge 1\right\}$ is a CONS in $H$. Therefore,
\begin{equation}\lim_{t\to T}\mbox{\rm tr}\left(\left(I-V_t^{*}V_
t\right)Q_T\right)=0\label{z1}\end{equation}
by
Lemma \ref{vt1} and the Dominated Convergence Theorem.
Combining (\ref{z0}) and (\ref{z1}) we obtain (\ref{c1}) and, consequently, (c).
Part (d) follows immediately from Lemma \ref{vt1a} and (\ref{covp}).
\end{proof}
\begin{proposition}\label{CL}
The conditional distribution of the process $(Z^x_t)$ in the space $H_2 = L^2(0,T;H)$ given
$Z^x_T$ is $N(\lambda , \overline{Q})$, where
\begin{equation}
\lambda (t) = S_t x + K_tQ^{-1/2}_T Z_T, \label{X1}
\end{equation}
\begin{equation}
\overline{Q} = \tilde Q - \mathcal{K} \mathcal {K}^*, \label{X2}
\end{equation}
where
$\tilde Q$ is the covariance operator of the process $(Z^x_t)$ in $H_2$, $\tilde{Q}: H\to H_2$,
$$[\tilde Q y] (t) = \int_0^t R(t,s)y(s)ds,\quad y\in H_2,$$
and
$$R(t,s)z = \int_0^s S_{t-r}QS^*_{s-r}z dr,\quad z\in H,\quad 0\le s\le t \le T,$$
and $\mathcal{K} : H\to H_1$ is defined in (\ref{dk}).
\end{proposition}
\begin{proof}
We use Theorem \ref{tcond} with $H_1= H$, $H_2= L^2(0,T;H)$, $X=Z^x_t$, $Y=(Z^x_t)$, $C_X=Q_T$, and
$C_Y=\tilde Q$. By the definition of the covariance $C_{XY}$,
$$\left\langle C_{XY}k, h \right\rangle _{L^2(0,T;H)} = \mathbf{E}
\left\langle Z^x_T,k\right\rangle\left\langle Z^x,h\right\rangle _{L^2(0,T,H)}, \quad k\in H_1,\, h\in H_2,$$
it is easy to compute $[C_{XY} k](t) = Q_tS^*_{T-t}k,\, t\in [0,T]$. Hence we have
$T^* = \overline{C_{XY}C^{-1/2}_X} = \mathcal{K}$ and $T:H_2 \to H_1,\,
Ty = \mathcal{K}^* y = \int_0^T K_t^*y(t) dt$. By Theorem \ref{tcond} we have that
$$\overline{Q} = C_Y - T^* T = \tilde Q -  \mathcal{K} \mathcal {K}^*,$$
 and
$$\lambda (t) = \mathbb{E} (Z_t^x|Z_T^x) = \mathbb{E} (S_tx + Z_t | Z^x_T) =\mathbb{E}(S_tx +\hat Z_t +
 K_tQ^{-1/2}_tZ_T|Z^x_T)$$
which yields $\lambda (t) = S_tx + K_tQ^{-1/2}_T Z_T$, because $\hat Z_t$ and $Z^x_T$ are stochastically independent,
hence (\ref{X1}) and (\ref{X2}) hold true.
\end{proof}
Recall that $\mu _T$ denotes the probability law of $Z_T$ on $H$.
\begin{proposition}\label{tp2}
There exists a Borel subspace $\mathcal M\subset H$ such that $\mu_T(\mathcal M)=1$ and for all $x\in H$ and $y\in
S_Tx+\mathcal M$ the $H$-valued Gaussian process
\begin{equation}\hat {Z}_t^{x,y}=Z_t^x- \mathcal{K}Q_T^{-1/2}\left(Z_T^x-y\right
),\label{p0}\end{equation}
is well defined with paths in $L^2(0,T;H)$ and
\begin{equation}
\hat {Z}_t^{x,y}=S_tx- \mathcal{K}Q_T^{-1/2}\left(S_Tx-y\right)+\hat {Z}_
t,\quad\mathbb P-a.s.\label{p1}
\end{equation}
\end{proposition}
\begin{proof}
By Lemma \ref{tmeas} we can choose a measurable linear space $\mathcal M$ such that $\mathcal{K}Q_T^{-1/2}$ is linear
on $\mathcal M$ and $\mu_T\left(\mathcal M\right)=1$. Therefore, $\mathcal{K}Q_T^{-1/2}\left(Z_T^x-y\right)$ is well
defined for any $y\in S_Tx+\mathcal M$ and (\ref{p1}) holds.
\end{proof}
\begin{theorem}\label{repr}
Let
$\Phi :L^2(0,T;H)\to\mathbb R$ be a
Borel mapping such that
\[\mathbb E\left|\Phi\left(Z^x\right)\right|<\infty .\]
Then
\begin{equation}\label{DOU}
\mathbb E\left(\left.\Phi\left(Z^x\right)\right|Z_T^x=y\right)=\mathbb E\Phi\left(\hat Z^{x,y}\right),\quad\mu_T^x-a.e.
\end{equation}
where the left-hand side of (\ref{DOU}) is defined as a function $g_\Phi = g_\Phi (y) \in L^1(H,\mu_T^x)$ such that
 $\mathbb{E}(\Phi (Z^x)|Z^x_T) = g_\Phi (Z^x_T)\ \mathbb{P}$-a.s.
\end{theorem}
\begin{proof}
We have to show that
$$\mathbb{E}(\Phi (Z^x)|Z^x_T) = \mathbb{E} (\Phi (\hat Z^{x,y}))|_{Z^x_T = y}\quad \mathbb{P}-a.s.$$
By Proposition \ref{CL} we have
\begin{equation}\label{conde1}
\mathbb{E}(\Phi (Z^x)|Z^x_T)= \int_H \Phi (z) N(\lambda , \overline{Q}) (dz) \quad \mathbb{P}-a.s.,
\end{equation}
where $\lambda$ and $\overline{Q}$ are defined by (\ref{X1}) and (\ref{X2}), respectively.
On the other hand, the covariance operator $\hat Q$ of the process $\hat Z^{x,y}_t$ in $H_2$ is by (\ref{p1}) the same as
the one of $\hat Z_t$. Since $Z_t= \hat Z_t + K_tQ_T^{-1/2}Z_T$ and the summands on the right-hand side are independent random variables, we obtain $\tilde Q = \hat Q +\mathcal{K} \mathcal {K}^*$,
that is, $\hat Q = \overline{Q}$. Also, we have
$$ \mathbb{E} \hat Z^{x,y}_t =S_tx - \mathcal{K}Q_T^{-1/2}(S_Tx - y),$$
and therefore
$$\mathbb{E} (\Phi (\hat Z^{x,y}))|_{Z^x_T = y} = \int_H \Phi (z) N(S_tx - K_t Q_T^{-1/2}(S_Tx -y), \overline{Q})(dz)|_{Z^x_T =y}
=\int_H \Phi (z) N(\lambda, \overline{Q})(dz)$$ $\mathbb{P}-a.s.$, which together with (\ref{conde1}) concludes the
proof.
\end{proof}
\begin{definition}\label{OUB}
Given $x,\ y\in H$ and an $H$-valued OU process $(Z_t^x)$, a process $(\hat Z_t^{x,y})$ satisfying (\ref{DOU}) is called an Ornstein-Uhlenbeck Bridge (connecting points $x$ at time $t=0$ and $y$ at time $t=T$). The probability law of the process $(\hat Z_t^{x,y})$ in the space $L^2(0,T;H)$ will be denoted by $\hat{\mu}^{x,y}$.
\end{definition}
Thus we have shown that the OU Bridge may be written in the form (\ref{p0}) or (\ref{p1}) and its probability law
$\hat{\mu}^{x,y}$ is $N(\gamma ,\overline Q)$ where $\gamma (t) = \mathbb{E}[\lambda (t)| Z^x_T =y] = S_tx - K_tQ^{-1/2}_T(S_Tx-y)$ $\mu _T^x - a.e.$
\par\noindent
The following Theorem has been proved in \cite{den} :
\begin{theorem}\label{conc}
Let $\mathcal E$ be a Banach space such that $\mu (\mathcal E)=1$. Then
$\hat{\mu}^{0,y}(\mathcal E)=1$ for $y\in\mathcal M$.
\end{theorem}
\section{SDE associated to the OU Bridge}
In the sequel we will need the following
\begin{hypothesis}\label{h1}
For any $t>0$
\begin{equation}\mbox{\rm im}\left(S_tQ^{1/2}\right)\subset\mbox{\rm im}\left
(Q_t^{1/2}\right).\label{im}\end{equation}
\end{hypothesis}
\begin{remark}\label{r1}
Condition (\ref{im}) is satisfied in some important cases.
\par\noindent
(a) If the process $(Z^x_t)$ is strong Feller
then $\mathrm{im}\left(S_t\right)\subset\mathrm{im}\left(Q_t^{1/2}\right
)$ and
therefore (\ref{im}) holds.
\par\noindent
(b) Let $H_Q=Q^{1/2}(H)$ be endowed with the norm $|x|_Q=\left|Q^{ -1/2}x\right|$, where $Q$ is assumed to be nondegenerate. Assume that $
S_tH_Q\subset H_Q$ for all $t\ge 0$ and $\left(S_t\right)$ restricted to $ H_Q$ is a $C_0$-semigroup. It was proved in \cite{acta} that in this case
$ S_t(H)\subset Q_t^{1/2}(H)$ for all $t>0$ and there exists $c>0$ such that
\[\left\|Q_t^{-1/2}S_tQ^{1/2}\right\|\le\frac c{\sqrt {t}},\quad
t>0.\]
Assume additionally that the process $(Z^x_t)$ admits a Gaussian
invariant measure $\nu$. Then, cf. \cite{acta},
$\left(S_t\right)$ is a $C_0$-semigroup on $H_Q$ if the transition semigroup
of the process $(Z^x_t)$ is analytic on $L^2(H,\nu )$, in particular this holds
for a symmetric Ornstein-Uhlenbeck process. Explicit
conditions for the analyticity and symmetry of the transition
semigroup of the process $(Z^x_t)$ in $L^2(H,\nu )$ may be found in \cite{acta}
and \cite{symm}.
\end{remark}
\begin{lemma}\label{lh1}
Assume that Hypothesis \ref{h1} holds. Then the function
\[t\to\left|Q_t^{-1/2}S_tQ^{1/2}h\right|,\]
is nonincreasing on $(0,\infty )$ for each $h\in H$.
\end{lemma}
\begin{proof}
By Lemma \ref{vt1}\ we have
\begin{equation}\left\|Q_{t+s}^{-1/2}S_tQ_s^{1/2}\right\|\le 1.\label{jan1}\end{equation}
By assumption the operator $Q_{t+s}^{-1/2}S_{t+s}Q^{1/2}$ is well defined and
bounded and $S_sQ^{1/2}h\in\mathrm{im}\left(Q_s^{1/2}\right)$. Therefore,
by (\ref{jan1})
\[\left|Q_{t+s}^{-1/2}S_{s+t}Q^{1/2}h\right|=\left|Q_{t+s}^{-1/2}
S_tQ_s^{1/2}Q_s^{-1/2}S_sQ^{1/2}h\right|\]
\[\le\left|Q_s^{-1/2}S_sQ^{1/2}h\right|,\]
and (b) follows.
\end{proof}
Let
\[Y_u=\int_u^TS_{T-s}Q^{1/2}dW_s,\quad u\le T.\]
Since the operator-valued function $t\to Q_t$ is continuous in the
weak operator topology and all the operators $Q_t$ are compact
for $t>0$, there exists a measurable choice of
eigenvectors $\left\{e_k(t):k\ge 1\right\}$ and eigenvalues  $\left
\{\lambda_k(t):k\ge 1\right\}$.
For each $n\ge 1$ we define a process
\[X^n_u=\sum_{k=1}^n\frac 1{\sqrt{\lambda_k(T-u)}}\left\langle Y_u,e_k(T-u)\right\rangle F_u^{*}e_k(T-u),\]
where $F_u=Q_{T-u}^{-1/2}S_{T-u}Q^{1/2}$.
\begin{lemma}\label{yu}
There exists a measurable stochastic process $\left(X_u\right)$ defined on $[0,T)$ such that for each $ a<T$
\begin{equation}\lim_{n\to\infty}\mathbb E\int_0^a\left|X_u^n-X_u\right
|^2du=0.\label{eps}\end{equation}
and for each $h\in H$ and $a<T$ the series
\begin{equation}\left\langle X_u,h\right\rangle =\sum_{k=1}^{\infty}\frac
1{\sqrt {\lambda_k(T-u)}}\left\langle Y_u,e_k(T-u)\right\rangle\left
\langle e_k(T-u),F_uh\right\rangle\label{fu}\end{equation}
converges in $L^2(0,a)$ in mean square.
Moreover, if $0\le u\le v<T$ then for all $h,k\in H$
\begin{equation}\mathbb E\left\langle X_u,h\right\rangle\left\langle
X_v,k\right\rangle =\left\langle F_uh,Q_{T-u}^{-1/2}Q_{T-v}^{1/2}
F_vk\right\rangle ,\label{covx}\end{equation}
where the operator $Q_{T-u}^{-1/2}Q_{T-v}^{1/2}$ is bounded.
\end{lemma}
\begin{proof}
For $u\le v\le T$
\begin{equation}
\mathbb E\left\langle Y_u,h\right\rangle\left\langle
Y_v,k\right\rangle =\left\langle Q_{T-v}h,k\right\rangle ,\quad h
,k\in H.\label{q}
\end{equation}
Therefore
\[\mathbb E\left\langle X_u^n-X_u^m,h\right\rangle^2=\sum_{j=m+1}^
n\frac 1{\lambda_k(T-u)}\mathbb E\left\langle Y_u,e_k(T-u)\right\rangle^
2\left\langle e_k(T-u),F_uh\right\rangle^2\]
\begin{equation}
=\sum_{j=m+1}^n\left\langle e_k(T-u),F_uh\right\rangle^
2\underset{n,m\to\infty}\longrightarrow 0,\label{koszi}
\end{equation}
hence the process
\[\left\langle X_u,h\right\rangle =\sum_{k=1}^{\infty}\frac 1{\sqrt {
\lambda_k(T-u)}}\left\langle Y_u,e_k(T-u)\right\rangle\left\langle
e_k(T-u),F_uh\right\rangle =\left\langle Q_{T-u}^{-1/2}Y_u,F_uh\right
\rangle\]
is well defined for each $h\in H$ and $u<T$. For $u,v$ such
that $0<u\le v<T$ we have
\begin{equation}\mathrm{im}\left(Q_{T-v}^{1/2}\right)\subset\mathrm{
im}\left(Q_{T-u}^{1/2}\right).\label{im1}\end{equation}
Let $P_n$ is an orthogonal projection on
$\mbox{\rm lin}\left\{e_k(T-v):k\le n\right\}$ and $F_u^n=P_nF_u$. Then $
Q_{T-u}^{-1/2}F^n_u$ is
bounded on $H$. Let
\[X_u^n=\left(Q_{T-u}^{-1/2}F_u^n\right)^{*}Y_u.\]
By (\ref{q})
\[\mathbb E\left\langle X_u^n,h\right\rangle\left\langle X_v^nk\right
\rangle =\left\langle Q_{T-v}Q_{T-u}^{-1/2}F_u^nh,Q_{T-v}^{-1/2}F_
v^nk\right\rangle\]
\[=\left\langle F_u^nh,Q_{T-u}^{-1/2}Q_{T-v}^{1/2}F_v^nk\right\rangle
.\]
By (\ref{im1}) the operator $Q_{T-u}^{-1/2}Q_{T-v}^{1/2}$ is bounded and
therefore
\[\mathbb E\left\langle Q_{T-u}^{-1/2}Y_u,F_uh\right\rangle\left\langle
Q_{T-v}^{-1/2}Y_v,F_vk\right\rangle =\lim_{n\to\infty}\mathbb E\left
\langle X_u^n,h\right\rangle\left\langle X_v^n,k\right\rangle\]
\[=\left\langle F_uh,Q_{T-u}^{-1/2}Q_{T-v}^{1/2}F_vk\right\rangle
.\]
It follows from (\ref{covx}) that
\[\mathbb E\left\langle X_u,h\right\rangle^2=\left|F_uh\right|^2,\]
and by Lemma \ref{lh1} we obtain for $u\le a$
\[\mathbb E\left\langle X_u^n,h\right\rangle^2\le\mathbb E\left\langle
X_u,h\right\rangle^2\le |h|^2\left\|F_{T-a}\right\|^2.\]
Then (\ref{koszi}) and the Dominated Convergence Theorem yield
\[\lim_{n,m\to\infty}\int_0^a\sup_{|h|\le 1}\mathbb E\left\langle
X_u^n-X_u^m,h\right\rangle^2du=0.\]
As a consequence we find that (\ref{eps}) holds for any
$a\in (0,T)$.
\end{proof}
By Lemma \ref{lh1}
a cylindrical process
\[I_t=\int_0^tF_u^{*}Q_{T-u}^{-1/2}Y_udu\]
is well defined, that is for any $h\in H$ the real-valued
process
\[\left\langle I_t,h\right\rangle =\int_0^t\left\langle Q_{T-u}^{-1/2}Y_u,F_uh\right\rangle du\]
is well defined for all $t<T$.

\begin{lemma}\label{wiener}
The cylindrical process
\[\zeta_t=W_t-\int_0^tF_u^{*}Q_{T-u}^{-1/2}Y_udu,\quad t\le T,\]
is a standard cylindrical Wiener process on $H$.
\end{lemma}
The proof of this Lemma is omitted; it is a word by word repetition of the proof of Lemma 4.7 in \cite{den} if we use
Lemmas \ref{lh1} and \ref{yu} above.
\begin{theorem}\label{weak1}
For all $t<T$
\begin{equation}
\mathbb E\int_0^t\left|S_{t-s}Q^{1/2}F_s^{*}Q_{T-
s}^{-1/2}S_{T-s}\hat Z_s\right|^2ds<\infty ,\label{pa}
\end{equation}
and
\begin{equation}
\hat {Z}_t=-\int_0^tS_{t-s}Q^{1/2}F_s^{*}Q_{T-s}^{
-1/2}S_{T-s}\hat {Z}_sds+\int_0^tS_{t-s}Q^{1/2}d\zeta_s,\quad\mathbb
P-a.s.\label{main}
\end{equation}
\end{theorem}
\begin{proof}
We will show first that the operator
$Q_{T-s}^{-1/2}S_{T-s}\hat {Q}_sS_{T-s}^{*}Q_{T-s}^{-1/2}$ is bounded. Let $
h,k\in H$. Then by
Proposition \ref{ou0}\  and (\ref{qtd}) we obtain
\[\left\langle S_{T-s}\hat Q_sS_{T-s}^{*}h,k\right\rangle =\left\langle
S_{T-s}Q_sS_{T-s}^{*}h,k\right\rangle -\left\langle S_{T-s}Q_s^{1
/2}V_s^{*}V_sQ_s^{1/2}S_{T-s}^{*}h,k\right\rangle\]
\[=\left\langle\left(Q_T-Q_{T-s}\right)h,k\right\rangle -\left\langle
Q_T^{-1/2}S_{T-s}Q_sS_{T-s}^{*}h,Q_T^{-1/2}S_{T-s}Q_sS_{T-s}^{*}k\right
\rangle\]
\[=\left\langle\left(Q_T-Q_{T-s}\right)h,k\right\rangle -\left\langle
Q_T^{-1/2}\left(Q_T-Q_{T-s}\right)h,Q_T^{-1/2}\left(Q_T-Q_{T-s}\right
)k\right\rangle\]
\[=\left\langle\left(Q_T-Q_{T-s}\right)h,k\right\rangle -\left\langle\left
(Q_T-Q_{T-s}\right)Q_T^{-1}\left(Q_T-Q_{T-s}\right)h,k\right\rangle\]
\[=\left\langle\left(Q_T-Q_{T-s}\right)h,k\right\rangle -\left\langle\left
(Q_T-Q_{T-s}\right)\left(I-Q_T^{-1}Q_{T-s}\right)h,k\right\rangle\]
\[=\left\langle\left(Q_{T-s}-Q_{T-s}Q_T^{-1}Q_{T-s}\right)h,k\right
\rangle =\left\langle Q_{T-s}^{1/2}\left(I-Q_{T-s}^{1/2}Q_T^{-1}Q_{
T-s}^{1/2}\right)Q_{T-s}^{1/2}h,k\right\rangle .\]
Since the operator $Q_{T-s}^{1/2}Q_T^{-1}Q_{T-s}^{1/2}$ is bounded for $
s<T$ we
find that the operator
\begin{equation}T_s=Q_{T-s}^{-1/2}S_{T-s}\hat {Q}_sS_{T-s}^{*}Q_{
T-s}^{-1/2}=I-Q_{T-s}^{1/2}Q_T^{-1}Q_{T-s}^{1/2}\label{ts}\end{equation}
is bounded as well . Therefore, for $s\le T-\epsilon$ Lemma
\ref{lh1} and (\ref{ts}) yield
\[\mathbb E\left|S_{t-s}Q^{1/2}F_s^{*}\left(Q_{T-s}^{-1/2}S_{T-s}
\hat Z_s\right)\right|^2\]
\[=\left\|S_{t-s}Q^{1/2}F_s^{*}\left(Q_{T-s}^{-1/2}S_{T-s}\hat Q_
s^{1/2}\right)\right\|_{HS}^2\le\left\|S_{t-s}Q^{1/2}\right\|_{HS}^
2\left\|F_s\right\|^2\left\|T_s^{1/2}\right\|\]
\[\le\left\|S_{t-s}Q^{1/2}\right\|_{HS}^2\left\|F_{T-\epsilon}\right
\|^2,\]
which completes the proof of (\ref{pa}). As a byproduct of
the argument given above we proved also that the
process $Q_{T-s}^{-1/2}S_{T-s}\hat {Z}_s$ is well defined for all $
s\le T$.
Now, we are ready to prove (\ref{main}). By Lemma
\ref{wiener}\ we have
\[\hat {Z}_t=Z_t-K_tQ_T^{-1/2}Z_T\]
\[=\int_0^tS_{t-s}Q^{1/2}d\zeta_s+\int_0^tS_{t-s}Q^{1/2}F_s^{*}Q_{
1-s}^{-1/2}Y_sds-K_tQ_T^{-1/2}Z_T,\]
and since
\[Y_s=Z_T-S_{T-s}Z_s=Z_T-S_{T-s}K_sQ_T^{-1/2}Z_T-S_{T-s}\hat {Z}_
s,\]
we find that
\[\hat {Z}_t=\int_0^tS_{t-s}Q^{1/2}d\zeta_s-\int_0^tS_{t-s}Q^{1/2}
F_s^{*}Q_{T-s}^{-1/2}S_{T-}{}_s\hat {Z}_sds\]
\[+\int_0^tS_{t-s}Q^{1/2}F_s^{*}Q_{T-s}^{-1/2}\left(Z_T-S_{T-s}K_
sQ_T^{-1/2}Z_T\right)ds-K_tQ_T^{-1/2}Z_T.\]
It remains to show that
\begin{equation}\int_0^tS_{t-s}Q^{1/2}F_s^{*}Q_{T-s}^{-1/2}\left(
Z_T-S_{T-s}K_sQ_T^{-1/2}Z_T\right)ds-K_tQ_T^{-1/2}Z_T=0.\label{00}\end{equation}
To this end note
first that
\begin{equation}K_tQ_T^{-1/2}Z_T=\left(\int_0^tS_{t-s}Q^{1/2}F_s^{
*}ds\right)Q_T^{-1/2}Z_T,\label{k1}\end{equation}
and
\[S_{T-t}K_tQ_T^{-1/2}Z_T=\left(\int_0^tS_{T-s}Q^{1/2}F_s^{*}ds\right
)Q_T^{-1/2}Z_T\]
\[=\left(Q_T-Q_{T-t}\right)Q_T^{-1}Z_T=Z_T-Q_{T-t}Q_T^{-1}Z_T,\]
and thereby
\begin{equation}Z_T-S_{T-t}K_tQ_T^{-1/2}Z_T=Q_{T-t}Q_T^{-1}Z_T.\label{k2}\end{equation}
Finally, (\ref{k2}) and the definition of $F_s^{*}$ give
\[\int_0^tS_{t-s}Q^{1/2}F_s^{*}Q_{T-s}^{-1/2}\left(Z_T-S_{T-s}K_s
Q_T^{-1/2}Z_T\right)ds\]
\[=\int_0^tS_{t-s}Q^{1/2}F_s^{*}Q_{T-s}^{-1/2}Q_{T-s}Q_T^{-1}Z_Td
s=\left(\int_0^tS_{t-s}Q^{1/2}F_s^{*}ds\right)Q_T^{-1/2}Z_T,\]
and (\ref{00}) follows from (\ref{k1}).
\end{proof}
We will consider now the general case of the bridge $\left(\hat Z_t^{x,y}\right)$ connecting points $x\in H$ and $y$.
We will impose the stronger condition (\ref{hh2}) which is now formulated as a separate hypothesis:
\begin{hypothesis}\label{h2}
For any $t\in (0,T]$
\[\mbox{\rm im}\left(Q_t^{1/2}\right)=\mbox{\rm im}\left(Q_T^{1/2}\right
).\]
\end{hypothesis}
For $y\in H_1 := \mbox{\rm im}(Q_T^{1/2})$ we define
\[Ny(t)=\int_0^tS_{t-s}Q^{1/2}F_s^{*}Q_{T-s}^{-1/2}yds,\quad t\le
T-\epsilon .\]
\begin{theorem}\label{vy}
Assume that Hypotheses \ref{h1} and \ref{h2} hold. Then
the following holds.
\par\noindent
(a) The
operator $N:H_1\to L^2\left(0,T-\epsilon ;H\right)$ is Hilbert-Schmidt.
\par\noindent
(b) For any $x\in H$ and $y\in\mathcal M$
\[\hat {Z}_t^{x,y}=S_tx-\int_0^tS_{t-s}Q^{1/2}F_s^{*}Q_{T-s}^{-1/
2}S_{T-s}\hat {Z}_s^{x,y}ds+\int_0^tS_{t-s}Q^{1/2}d\zeta_s\]
\begin{equation}+\int_0^tS_{t-s}Q^{1/2}F_s^{*}Q_{T-s}^{-1/2}yds.\label{b2}\end{equation}
\end{theorem}
\begin{proof}
Recall that by Lemma \ref{ania} (b) we have $\mathcal{K} =\mathcal{L} B$,
hence for $z\in\mathcal M$
\begin{equation}K_tQ_T^{-1/2}z=\int_0^tS_{t-s}Q^{1/2}B_sQ_T^{-1/2}
zds.\label{ky}\end{equation}
Next, for $s\le T-\epsilon$
\[\sup_{s\le T-\epsilon}\left\|Q_{T-s}^{-1/2}Q_T^{1/2}\right\|=\left
\|Q_{\epsilon}^{-1/2}Q_T^{1/2}\right\|<\infty ,\]
and invoking Lemma \ref{lh1} we find that
\[\left\|NQ_T^{1/2}\right\|_{HS}^2\le\int_0^{T-\epsilon}\left\|\int_
0^tS_{t-s}Q^{1/2}F_s^{*}Q_{T-s}^{-1/2}Q_T^{1/2}ds\right\|_{HS}^2d
t\]
\[\le\left(\int_0^T\left\|S_sQ^{1/2}\right\|^2_{HS}ds\right)\left
(\int_0^{T-\epsilon}\left\|F_s^{*}Q_{T-s}^{-1/2}Q_T^{1/2}\right\|^
2ds\right)\]
\[\le\left\|F_{T-\epsilon}\right\|^2\left\|Q_{\epsilon}^{-1/2}Q_T^{
1/2}\right\|^2\left(\int_0^T\left\|S_sQ^{1/2}\right\|^2_{HS}ds\right
)<\infty .\]
Therefore, the measurable function
\[y\to\int_0^tS_{t-s}Q^{1/2}F_s^{*}Q_{T-s}^{-1/2}yds,\]
is well defined.
We are ready now for the proof of (\ref{b2}). Let
$x,y\in\mbox{\rm im}\left(Q_T^{1/2}\right)$. Then Hypothesis \ref{h2} yields
$S_Tx\in\mbox{\rm im}\left(Q_T^{1/2}\right)$, hence $y\in\mathcal
M$. By (\ref{p1}) we have
\[\hat {Z}_t^{x,y}=\hat {Z}_t+S_tx-K_tQ_T^{-1/2}\left(S_Tx-y\right
),\]
and Theorem \ref{weak1}\ yields
\[\hat {Z}_t^{x,y}=S_tx-K_tQ_T^{-1/2}S_Tx+K_tQ_T^{-1/2}y\]
\[-\int_0^tS_{t-s}Q^{1/2}F_s^{*}Q_{T-s}^{-1/2}S_{T-s}\hat {Z}_sds
+\int_0^tS_{t-s}Q^{1/2}d\zeta_s\]
\[=S_tx-K_tQ_T^{-1/2}S_Tx+K_tQ_T^{-1/2}y\]
\[-\int_0^tS_{t-s}Q^{1/2}F_s^{*}Q_{T-s}^{-1/2}S_{T-s}\left(\hat Z_
s^{x,y}-S_sx+K_sQ_T^{-1/2}S_Tx-K_sQ_T^{-1/2}y\right)ds+\int_0^tS_{
t-s}Q^{1/2}d\zeta_s\]
\[=-K_tQ_T^{-1/2}S_Tx+\int_0^tS_{t-s}Q^{1/2}F_s^{*}Q_{T-s}^{-1/2}
S_{T-s}\left(S_s-K_sQ_T^{-1/2}S_T\right)xds\]
\[+K_tQ_T^{-1/2}y+\int_0^tS_{t-s}Q^{1/2}F_s^{*}Q_{T-s}^{-1/2}S_{T
-s}K_sQ_T^{-1/2}yds\]
\[+S_tx-\int_0^tS_{t-s}Q^{1/2}F_s^{*}Q_{T-s}^{-1/2}S_{T-s}\hat {Z}_
s^{x,y}ds+\int_0^tS_{t-s}Q^{1/2}d\zeta_s\]
\begin{equation}=:H_tx+G_ty+S_tx-\int_0^tS_{t-s}Q^{1/2}F_s^{*}Q_{T
-s}^{-1/2}S_{T-s}\hat {Z}_s^{x,y}ds+\int_0^tS_{t-s}Q^{1/2}d\zeta_
s.\label{hg}\end{equation}
We will show first that
\begin{equation}G_ty=\int_0^tS_{t-s}Q^{1/2}F_s^{*}Q_{T-s}^{-1/2}y
ds.\label{gt}\end{equation}
For $y\in\mbox{\rm im}\left(Q_T^{1/2}\right)$
\[S_{T-t}K_ty=\int_0^tS_{t-s}QS_{T-s}^{*}Q_T^{-1/2}yds=\int_0^tS_{
T-s}QS_{T-s}^{*}Q_T^{-1/2}yds\]
\begin{equation}=\left(Q_T-Q_{T-t}\right)Q_T^{-1/2}y,\label{s1}\end{equation}
and therefore
\[F_s^{*}Q_{T-s}^{-1/2}S_{T-s}K_sy=F_s^{*}Q_{T-s}^{-1/2}Q_T^{1/2}
y-F_s^{*}Q_{T-s}^{1/2}Q_T^{-1/2}y\]
\[=F_s^{*}Q_{T-s}^{-1/2}Q_T^{1/2}y-Q^{1/2}S_{T-s}^{*}Q_T^{-1/2}y.\]
Hence, taking Lemma \ref{ania} (b) into account we find that
\[G_ty=K_tQ_T^{-1/2}y+\int_0^tS_{t-s}Q^{1/2}F_s^{*}Q_{T-s}^{-1/2}
S_{T-s}K_sQ_T^{-1/2}yds\]
\[=K_tQ_T^{-1/2}y+\int_0^tS_{t-s}Q^{1/2}F_s^{*}Q_{T-s}^{-1/2}yds-
K_tQ_T^{-1/2}y,\] and (\ref{gt}) follows. Next, we claim that for $x\in\mbox{\rm im}\left(Q_T^{1/2}\right)$
\begin{equation}H_tx=0.\label{ht}\end{equation}
Indeed, using (\ref{s1}) we obtain
\[H_tx=-K_tQ_T^{-1/2}S_Tx+\int_0^tS_{t-s}Q^{1/2}F_s^{*}Q_{T-s}^{-
1/2}S_{T-s}\left(S_s-K_sQ_T^{-1/2}S_T\right)xds\]
\[=-K_tQ_T^{-1/2}x+\int_0^tS_{t-s}Q^{1/2}F_s^{*}Q_{T-s}^{-1/2}S_T
xds-\int_0^tS_{t-s}Q^{1/2}F_s^{*}Q_{T-s}^{-1/2}S_{T-s}K_sQ_T^{-1/
2}S_Txds\]
\[=-K_tQ_T^{-1/2}x+\int_0^tS_{t-s}Q^{1/2}F_s^{*}Q_{T-s}^{-1/2}S_T
xds\]
\[-\int_0^tS_{t-s}Q^{1/2}F_s^{*}Q_{T-s}^{-1/2}\left(Q_T-Q_{T-t}\right
)Q_T^{-1/2}S_Txds=0,\]
which yields (\ref{ht}) for $x\in\mbox{\rm im}\left(Q_T^{1/2}\right
)$ and therefore for
all $x\in H$.
Finally, combining (\ref{hg}),
(\ref{gt}) and (\ref{ht}) we obtain (\ref{b2}).
\end{proof}
\begin{corollary}\label{col1}
Assume Hypotheses \ref{h1}\ and \ref{h2}. Then for each
$t<T$, and $h\in\mbox{\rm dom}\left(A^{*}\right)$ and all $x\in H$ and $
y\in\mathcal M$
\[\left\langle\hat Z_t^{x,y},h\right\rangle =\left\langle x,h\right
\rangle +\int_0^t\left\langle\hat Z_s^{x,y},A^{*}h\right\rangle d
s-\int_0^t\left\langle F_s^{*}Q_{T-s}^{-1/2}S_{T-s}\hat Z_s^{x,y}
,Q^{1/2}h\right\rangle ds\]
\[+\int_0^t\left\langle F_s^{*}Q_{T-s}^{-1/2}y,Q^{1/2}h\right\rangle
ds+\left\langle\zeta_t,Q^{1/2}h\right\rangle .\]
\end{corollary}
\begin{proof}
On any interval $\left[0,T_0\right]$ with $T_0<T$ and for any $y\in\mathcal M$ the functions
\[s\to Q^{1/2}F_s^*Q_{T-s}^{-1/2}S_{T-s}\hat{Z}^{x,y}_s\quad\mathrm{and}\quad s\to Q^{1/2}F_s^*Q_{T-s}^{-1/2}y\]
are $\mathbb P$-a.s. Bochner integrable by Theorem \ref{vy} and therefore standard results about the equivalence of weak and strong solutions of
deterministic and stochastic evolution equations can be applied to prove the corollary, see for example \cite{ball} for deterministic
and \cite{ania-mild}, \cite{on} for stochastic versions.
\end{proof}
\section{Applications to Semilinear Equations}
In this Section, transition densities and Markov semigroups defined by semilinear stochastic equations are studied
using the OU Bridge.
Throughout the Section we assume (beside (\ref{hou})) that the OU process $(Z^x_t)$ is strongly Feller, that is,
the condition
\begin{equation}\label{SFP}
\mathrm{im}(S_t) \subset \mathrm{im} (Q_t^{1/2}),\quad t\in (0,T),
\end{equation}
is satisfied. Note that (\ref{SFP}) trivially implies the preceding Hypotheses \ref{h1} and \ref{h2} (or (\ref{hh2})). Let
$(\mathcal{P}, \| .\| _{var})$ denote the space of probability measures on the Borel sets of $H$ endowed with the metric of total variation.
We start from a simple proposition where some continuity properties of the OU Bridge are
given.
\begin{proposition}\label{cont}
(a) For each $t\in (0,T), y\in \mathcal{M} $, where $\mathcal M$ has been defined in Proposition \ref{tp2}, the mappings
\begin{equation}\label{co1}
x\mapsto \hat Z^{x,y}_t(\omega ), \quad H \to H,
\end{equation}
\begin{equation}\label{co2}
x\mapsto \hat Z^{x,y}(\omega ), \quad H \to L^2(0,T;H),
\end{equation}
are continuous for $\mathbb{P}$-almost all $\omega \in \Omega $, and the mapping
\begin{equation}\label{co3}
x \mapsto \hat \mu ^{x,y}_t,\quad H \to (\mathcal{P}, \| . \|_{var}),
\end{equation}
is continuous.
\par\noindent
(b) If, moreover, for each $t\in(0,T)$ we have $\overline{K_tQ^{-1/2}_T} \in \mathcal{L}(\hat H, H)$, where
$\hat H$ is a separable Banach space continuously embedded into $H$, then
the mapping $y\mapsto \hat Z^{x,y}_t(\omega )$ is  $\hat H \to H$ $\mathbb{P}$- a.s. continuous. Similarly, if
\begin{equation}
\overline{\mathcal{K} Q_T^{-1/2}} \in \mathcal{L}(\hat H, L^2(0,T;H))
\end{equation}
then $\mathcal{M} \supset \hat H$ and the mapping $y \mapsto \hat Z^{x,y}(\omega )$ is $\mathbb{P}$-a.s. $\hat H\to L^2(0,T;H))$ continuous.
\end{proposition}
\begin{proof}
(a) By (\ref{SFP}) we have that $S_Tx \in \mathrm{im} (Q_T^{1/2})$ for each $x\in H$ and hence $S_Tx\in \mathcal{M}$ by
construction of $\mathcal{M}$, hence $y\in \mathcal{M}$. Furthermore, (\ref{SFP}) implies that the mappings
$\mathcal{K}Q_T^{-1/2}S_T$ and $K_tQ^{-1/2}_TS_T,\, t\in (0,T]$, are in $\mathcal{L}(H,L^2(0,T;H))$ and
$\mathcal{L}(H)$, respectively, and (\ref{co1}) and (\ref{co2}) follow by (\ref{p1}). To show (\ref{co3}) we recall
Proposition \ref{ou0} and Lemma \ref{vt1a} , by which we have $\mathrm{im}(\hat Q_t^{1/2}) = \mathrm{im} (Q_t^{1/2})$.
Hence the measures $(\hat \mu ^{x,y}_T), x\in H,$ are equivalent and
\begin{equation}\label{psi}
\psi^y(t,x,z)=\frac {d\mu^{x,y}_t}{d\mu^{0,y}_t}(z)=\exp\left(-\frac
12\left|Q_t^{-1/2}S_tx\right|^2+\frac 12\left|Q_T^{-1/2}S_Tx\right
|^2+\left\langle Q_t^{-1/2}z,Q_t^{-1/2}S_tx\right\rangle\right).
\end{equation}
Indeed, by the Cameron-Martin formula we have
\[\psi^y(t,x,z)=\exp\left(-\frac 12\left|\hat Q_t^{-1/2}m\right|^
2+\left\langle\hat Q_t^{-1/2}z,\hat Q_t^{-1/2}m\right\rangle\right
),\]
where $m=Q_t^{1/2}\left(I-V_t^{*}V_t\right)Q_t^{-1/2}S_tx$. Then using
(\ref{covp}) we get (\ref{psi}) and the assertion easily follows.
\par\noindent
The proof of part (b) is completely analogous.
\end{proof}
\begin{remark}\label{psi1}
(a) The equivalent form of the density (\ref{psi}) is
\[\psi^y(t,x,z)=\exp\left(-\frac 12\left|\left(I-V_t^{*}V_t\right
)^{1/2}Q_t^{-1/2}S_tx\right|^2+\left\langle Q_t^{-1/2}z,Q_t^{-1/2}
S_tx\right\rangle\right).\]
(b) Note that the OU Bridge $(\hat Z^{x,y}_t)$ satisfies the SDE (\ref{b2}) which defines
an (inhomogeneous) Markov process on the interval $(0,T)$. By (\ref{co3}) this process is strongly Feller.
\end{remark}
\vskip1cm
Now consider a stochastic semilinear evolution equation of the form
\begin{equation}\label{SEM}
dX_t = AX_tdt + F(X_t)dt + \sqrt{Q}dW_t,\quad X_0=x\in H
\end{equation}
where $A$, $W_t$ and $Q$ are as before and $F: H\to H$ is a nonlinear continuous mapping. Suppose that
$\mathrm{im}(F) \subset \mathrm{im} (Q^{1/2})$ and set $G:= Q^{-1/2}F$.
\begin{hypothesis}\label{GIR}
The mapping $G: H\to H$ is bounded and continuous.
\end{hypothesis}
Now we formulate technical assumptions on the linear part of the equation. For simplicity of presentation, it is stated in the form that is verifiable in
examples and includes all assumptions made previously in the paper.
\begin{hypothesis}\label{h4}
Assume either
\par\noindent
(i) $\mathrm{dim} H <\infty $ or
\par\noindent
(ii) There exist $\alpha \in (0,1)$ and $\beta < \frac{1+\alpha }{2}$ such that
\[ \int_0^{T_0} t^{-\alpha} \| S_t Q^{1/2}\| ^2_{HS} dt <\infty\quad and \]
\[ \| Q_t^{-1/2}S_t\| \le \dfrac{c}{t^{\beta}},\quad t\in (0,T_0), \]
for some $c>0$ and $T_0>0$.
\end{hypothesis}
Conditions from (ii) are often used in the theory of stochastic equations and have been widely studied
(cf.\cite{dz1} or \cite{den}, see also the Example below). Note that Hypothesis \ref{h4} (ii) implies all previous assumptions made in the paper on the linear part of the equation (\ref{SEM}) (i.e., all except for Hypothesis \ref{GIR}).
\par\noindent
It is well known (see e.g. \cite{on} ) that under Hypotheses \ref{GIR} and \ref{h4} the equation (\ref{SEM}) defines an $H$-valued
Markov process induced by the mild formula
\begin{equation}\label{MF}
X_t = S_tx + \int_0^t S_{t-r}F(X_r)dr + \int_0^t S_{t-r}\sqrt{Q}d\widetilde W_r, \quad t\ge 0,
\end{equation}
where $\widetilde W_t$ is a standard cylindrical Wiener process on $H$ defined on a suitable probability space.
\par\noindent
Finally, we assume that the OU process defined by the linear equation (\ref{01}) has an invariant measure $\nu$ that will be used as a reference measure. This is equivalent to the condition
\begin{equation}\label{ea1}
\sup_{t>0} tr (Q_t) < \infty.
\end{equation}
If (\ref{ea1}) holds then $\nu$ is a centered Gaussian measure with the covariance operator
\[Q_{\infty}=\int_0^{\infty}S_tQS_t^*dt.\]
Moreover, it has been shown in \cite{fock} that $S_tQ_{\infty}^{1/2}(H)\subset Q_{\infty}^{1/2}(H)$ and the family of operators
\[S_0(t)=Q_{\infty}^{-1/2}S_tQ_{\infty}^{1/2},\quad t\ge 0,\]
defines a $C_0$-semigroup of contractions on $H$. Moreover, if part (ii) of Hypothesis \ref{h4} holds then $\left\|S_0(t)\right\|<1$ for all $t>0$.
\par\noindent
 Denote by $(P_t)$ the transition Markov semigroup defined by the equation (\ref{SEM}) and set
$$P(t,x, \Gamma )= P_t1_\Gamma (x),\, x\in H,\, t>0$$
and $\Gamma$ Borel sets in $H$, and
$$d(t,x,y) = \frac{P(t,x,dy)}{\nu (dy)}.$$
It is standard to see that the density $d$ exists, because Girsanov Theorem may be used to show the equivalence of measures $P(t,x,dy) \sim \mu
^x_t$, and $\mu ^x_t \sim \nu$ by (\ref{SFP}) (see e.g. \cite{den}).
\begin{theorem}\label{dens}
Let Hypotheses \ref{GIR}, \ref{h4} and (\ref{ea1}) be satisfied and let $T>0$ be fixed.
Then for $\nu$-almost all $y\in H$ the mapping $x \mapsto d(T,x,y)$ is
continuous on $H$.
\end{theorem}
\begin{theorem}\label{lp}
 Let Hypotheses \ref{GIR}, \ref{h4} and (\ref{ea1}) be satisfied. Then for $p>1$, $T>0$, we have
\[ P_T(L^p(H,\nu )) \subset \mathcal{C} (H), \]
that is, the semigroup $(P_t)$ maps the space $L^p(H,\nu )$ into the space of continuous functions on $H$.
\end{theorem}
For $p,q>1$ we introduce the notation
\[\left\|P_t\right\|_{p,q}=\left(\int_H\left(\int_Hd^{p'}(t,x,y)\nu(dy)\right)^{q/p'}\nu(dx)\right)^{1/q},\]
where $p'=\frac{p}{p-1}$. Note that $\left\|P_t\right\|_{2,2}$ is a Hilbert-Schmidt norm of $P_t$. Moreover, if $\left\|P_t\right\|_{p,q}<\infty$
then the operator $P_t:L^p(H,\nu)\to L^q(H,\nu)$ is compact. Under assumptions more general than ours necessary and sufficient conditions were given
in \cite{ania} for boundedness of the operator $P_t:L^p(H,\nu)\to L^q(H,\nu)$. In the theorem below we use different arguments based on the formula
for transition densities to show that a stronger property holds: $\left\|P_t\right\|_{p,q}<\infty$.
\begin{theorem}\label{HSS}
Let Hypotheses \ref{GIR}, \ref{h4} and (\ref{ea1}) be satisfied. Then for any fixed $T>0$ and $q>0$ satisfying
\[q<1+\frac{p-1}{\left\|S_0(T)\right\|^2}\]
we have $\left\|P_T\right\|_{p,q}<\infty$. In particular, the operator $P_T:L^p(H,\nu)\to L^q(H,\nu)$ is $q$-summing and $P_T$ is
Hilbert-Schmidt in the space $L^2(H,\nu )$.
\end{theorem}
By the above mentioned equivalence of probabilities we may write
\begin{equation}\label{y1}
d(T,x,y) =\frac{P(T,x,dy)}{\mu^x_T(dy)}\cdot \frac{\mu^x_T(dy)}{\mu^0_T(dy)}\cdot \frac{\mu^0_T(dy)}{\nu(dy)}
\end{equation}
\begin{equation}
=: h(T,x,y)\cdot g(T,x,y)\cdot k(T,y),
\end{equation}
where $k$ does not depend on $x$, $g$ is given by the Cameron-Martin formula
\begin{equation}\label{y1a}
g(T,x,y) = \exp \{\left\langle x, \overline{S^*_TQ^{-1/2}_T} Q^{-1/2}_T y\right\rangle - \frac{1}{2}|Q^{-1/2}_TS_Tx|^2\}
\end{equation}
for $\nu$-almost all $y\in H$, and $h$ may be expressed by means of the OU Bridge $(\hat Z^{x,y}_t)$,
\begin{equation}\label{y2}
h(T,x,y) = \mathbb{E} \exp \{ \rho (\hat Z^{x,y} ) - \int_0^T \left\langle G(\hat Z^{x,y}_s), B_1(s) \hat Z_s +
B_2(s) x - B_3(s)y \right\rangle ds \}
\end{equation}
(cf.\cite{den}, Theorem 5.2), where
$$\rho (\hat Z^{x,y}) = \int_0^T \left\langle  G(\hat Z^{x,y}_s), dW_s\right\rangle
- \frac{1}{2}\int_0^T |G(\hat Z^{x,y}_s)|^2ds $$ and
$$ B_1(s) = (Q^{-1/2}_{T-s}S_{T-s}Q^{1/2})^* Q^{-1/2}_{T-s}S_{T-s},$$
$$B_2(s) = (Q^{-1/2}_TS_{T-s}Q^{1/2})^*Q_T^{-1/2}S_T, $$
$$B_3(s)y =(Q^{-1/2}_TS_{T-s}Q^{1/2})^*Q_T^{-1/2}y,\quad y\in \mathrm{im}\left(Q_T^{1/2}\right).$$
From Lemma \ref{ania} it follows that
\begin{equation}\label{y3}
\int_0^T |B_2(s)x|^2 ds = |Q^{-1/2}_TS_Tx|^2,\quad x\in H,
\end{equation}
and by \cite{den}, Proposition 4.9, we have that
\begin{equation}\label{y4}
\mathbb{E} \int_0^T |B_1(s)\hat Z_t| ds <\infty
\end{equation}
and
\begin{equation}\label{y5}
\int_0^T |B_3(s)y|ds <\infty
\end{equation}
for $\nu$- almost all $y \in \mathcal{M}$ (with no loss of generality we may assume that (\ref{y5}) holds for all
$y\in\mathcal{M},\, \nu (\mathcal{M}) =1$). The proofs of Theorems \ref{dens}, \ref{lp} and \ref{HSS} are based on the
following technical lemma:
\begin{lemma}\label{est}
Given $T>0$ and $q\in [0,\infty)$, there exists a constant $k_q>0$ such that
\begin{equation}\label{y6}
\begin{aligned}
h_q(T,x,y):&= \mathbb{E}\exp\{ q(\rho (\hat Z^{x,y} ) - \int_0^T \left\langle G(\hat Z^{x,y}_s), B_1(s) \hat Z_s + B_2(s) x - B_3(s)y \right\rangle
ds )\}\\
 &\le k_q \exp \{k_q (|x| + \int_0^T |B_3(s)y|ds )\}
 \end{aligned}
\end{equation}
for all $x\in H$ and $y\in \mathcal{M}$, in particular,
\[ h(t,x,y) \le k_1 \exp \{k_1 (|x| + \int_0^T |B_3(s)y|ds)\}. \]
\end{lemma}
\begin{proof}
By the Cauchy inequality we have
\begin{equation}\label{y7}
h_q(T,x,y) \le (\mathbb{E} \exp \{2q\rho (\hat Z^{x,y}\})^{1/2}
\end{equation}
$$\times (\mathbb{E} \exp \{2q(\int_0^T |\left\langle G(\hat Z^{x,y}_s), B_1(s) \hat Z_s +
B_2(s) x - B_3(s)y \right\rangle |ds)\} )^{1/2} $$
and since the process $s \mapsto G(\hat Z^{x,y}_s)$ is bounded the first expectation on the right-hand side
of (\ref{y7}) is bounded (uniformly w.r.t. $x$ and $y$). By (\ref{y3}) and (\ref{y5}) we thus have
\begin{equation}
h_q(T,x,y) \le C_q( \mathbb{E} \exp \{ C_q \int_0^T (|B_1(s)\hat Z_s| + |B_2(s)x| + |B_3(s)y|)ds\})^{1/2}
\end{equation}
$$ \le \tilde C_q \exp \{ \tilde C_q (|Q_t^{-1/2}S_Tx| + \int_0^T |B_3(s)y|ds)\} (\mathbb{E} \exp \{ \tilde C_q \int_0^T |B_1(s) \hat Z_s|ds \} )^{1/2} $$
for some $C_q,\, \tilde C_q$, and (\ref{y6}) follows by (\ref{y4}) and the Fernique inequality.
\end{proof}
\noindent{\it Proof of Theorem \ref{dens}.}
Without loss of generality (dropping, if necessary, a set of $\nu$-measure zero) we may suppose that $g(T,x,y)$ and $k(T,y)$ are defined for all $y\in \mathcal{M}$. By (\ref{y1a}) we have that the mapping $x\mapsto g(T,x,y)k(T,y)$ is continuous, so we only have to prove continuity of the mapping $x\mapsto h(T,x,y),\, y\in \mathcal{M},\, T>0$. Let $x_n \to x_0$ in $H$. First we show (possibly, for a subsequence) that
\begin{equation}\label{y8}
\lim _{n\to \infty} \exp \{ \rho (\hat Z^{x_n,y} ) - \int_0^T \left\langle G(\hat Z^{x_n,y}_s), B_1(s) \hat Z_s +
B_2(s) x_n - B_3(s)y \right\rangle ds\}
\end{equation}
$$ = \exp \{ \rho (\hat Z^{x_0,y} ) - \int_0^T \left\langle G(\hat Z^{x_0,y}_s), B_1(s) \hat Z_s +
B_2(s) x_0 - B_3(s)y \right\rangle ds\} $$
$\mathbb{P}$-a.s. We have
\begin{equation}\label{y9}
\begin{aligned}
&\int_0^T \left|\left\langle G(\hat Z^{x_n,y}_s), B_1(s) \hat Z_s + B_2(s) x_n - B_3(s)y \right\rangle  -  \int_0^T \left\langle G(\hat
Z^{x_0,y}_s),
B_1(s) \hat Z_s + B_2(s) x_0 - B_3(s)y \right\rangle\right| ds\\
&\le \int_0^T |G(\hat Z^{x_n,y}_s) - G(\hat Z^{x_0,y}_s)| (|B_1(s)\hat Z_s| + |B_2(s)x_0| + |B_3(s)y|)ds\\
&+ \int_0^T \left|G(\hat Z^{x_0,y}_s)|\cdot |B_2(s)(x_n-x_0)\right|ds,
\end{aligned}
\end{equation}
 which tends to zero by continuity and boundedness of $G$, (\ref{y3}) and Dominated
Convergence Theorem. Also, we have
\[ \mathbb{E} |\rho (\hat Z^{x_n,y}) - \rho (\hat Z^{x_0,y})| \le  C \left(\left(\mathbb{E} \int_0^T |G(\hat Z^{x_n,y}_s)
- G(\hat Z^{x_0,y}_s)|^2 ds\right)^{1/2}\right.\]
\[\left. + \mathbb{E} \int_0^T |G(\hat Z^{x_n,y}_s)  - G(\hat Z^{x_0,y}_s)|^2 ds \right),\]
which again tends to zero by Dominated Convergence Theorem, so there is a subsequence converging $\mathbb{P}$-a.s. Taking into account (\ref{y9}) we obtain (\ref{y8}). By (\ref{y6}) (used, for instance, with $q=2$) the random variables on the left-hand side of (\ref{y8}) are integrable uniformly in $n$, hence the convergence in (\ref{y8}) holds also in the space $L^1(\Omega )$ and, consequently, we obtain $h(T,x_n,y) \to h(T,x_0,y)$. Since we may choose a subsequence with this property from an arbitrary sequence $x_n \to x_0$, the convergence takes place for the whole sequence.
\par\medskip\noindent
{\it Proof of Theorem \ref{lp}.}
Let $T>0,\, \phi \in L^p(H, \nu)$ and $x_n \to x_0$ in $H$. Then
$$ |P_T\phi (x_n) - P_T \phi (x_0) | \le \int_H |\phi (y)| |d(T,x_n)-d(T,x_0,y)| \nu(dy) $$
$$\le (\int_H|\phi |^pd\nu )^{1/p} (\int_H |d(T,x_n, y) - d(T,x_0,y)|^{p'} \nu (dy))^{1/p'},$$
so by Theorem \ref{dens} it suffices to show that
\begin{equation}\label{y10}
\int_H (d(T,x_n,y))^q \nu (dy) < c_q,\quad q\in (1,\infty),
\end{equation}
where $c_q$ does not depend on $n$. The same property (uniform boundedness in arbitrary $L^q(H,\nu )$) has been shown for Gaussian densities $g(T,x_n,\cdot )$ and $k(T,\cdot )$ in \cite{reg}, so we only have to show
(\ref{y10}) where $d(T,x_n,y)$ is replaced by $h(T,x_n,y)$. However, by Lemma \ref{est} and H\" older inequality we have
\begin{equation}\label{y11}
 \int_H (h(T,x_n,y))^q \nu (dy) \le \int_H h_q(T,x_n,y) \nu (dy)
\end{equation}
$$\le k_q \exp \{k_q |x_n|\} \int_H \exp \{\int_0^T |B_3(s)y|ds\}\nu (dy) <c_q$$
where $c_q$ does not depend on $n$, since the sequence $x_n$ is obviously bounded and
$$ \int_H \exp \{\int_0^T |B_3(s)y|ds\}\nu (dy) <\infty $$
by (\ref{y5}), (\ref{SFP}) and the Fernique inequality.
\par\medskip\noindent
{\it Proof of Theorem \ref{HSS}.} We can rewrite (\ref{y1}) in the form
\[d(T,x,Y)=h(T,x,y)H(T,x,y),\]
where
\[H(T,x,y)=\frac{\mu_T^x(dy)}{\nu(dy)}.\]
Invoking the H\"older inequality we obtain
\begin{equation}\label{1015}
\begin{aligned}
\left\|P_T\phi\right\|_{p,q}^q&=\int_H\left(\int_H hH\phi\nu(dy)\right)^q\nu(dx)\\
&\le \int_H\left(\left(\int_Hh^{p'}H^{p'}\nu(dy)\right)^{1/p'}\left(\int_H|\phi |^p\nu(dy)\right)^{1/p}\right)^q\nu(dx)\\
&=\|\phi\|_p^q\int_H\left(\int_Hh^{p'}H^{p'}\nu(dy)\right)^{q/p'}\nu(dx).
\end{aligned}
\end{equation}
It remains to show that
\begin{equation}\label{1029}
K=\int_H\left(\int_Hh^{p'}H^{p'}\nu(dy)\right)^{q/p'}\nu(dx)<\infty .
\end{equation}
Indeed, using successively the H\"older equality we obtain for any $r>1$
\begin{equation}\label{637}
\begin{aligned}
K&\le \int_H\left(\int_Hh^{p'r'}\nu(dy)\right)^{q/p'r'}\left(\int_HH^{p'r}\nu(dy)\right)^{q/p'r}\nu(dx)\\
&\le \left(\int_H\left(\int_Hh^{p'r'}\nu(dy)\right)^{q/p'}\nu(dx)\right)^{1/r'}\left(\int_H\left(H^{p'r}\nu(dy)\right)^{q/p'}\nu(dx)\right)^{1/r}.
\end{aligned}
\end{equation}
It was shown in \cite{reg} that
\begin{equation}\label{655}
\int_H\left(H^{a'}\nu(dy)\right)^{b/a'}\nu(dx)<\infty ,
\end{equation}
for any $a,b\ge 1$, such that
\begin{equation}\label{813}
b\le 1+\frac{a-1}{\left\|S_0(T)\right\|^2}.
\end{equation}
Putting
\[a=\frac{p'r}{p'r-1}\quad\mathrm{and}\quad b=qr,\]
we find that there exists $r>1$ such that (\ref{813}) holds. Therefore, for such an $r$
\begin{equation}\label{816}
\int_H\left(H^{p'r}\nu(dy)\right)^{q/p'}\nu(dx)=\int_H\left(H^{a'}\nu(dy)\right)^{b/a'}\nu(dx)<\infty .
\end{equation}
 Next, we need to show that
 \begin{equation}\label{817}
 \int_H\left(\int_Hh^{p'r'}\nu(dy)\right)^{q/p'}\nu(dx)<\infty .
 \end{equation}
 To prove (\ref{817}) we note that if $\frac{q}{p'}\ge 1$ then
 \[\int_H\left(\int_Hh^{p'r'}\nu(dy)\right)^{q/p'}\nu(dx)\le \int_H\int_Hh^{r'q}\nu(dy)\nu(dx)\]
However, using Lemma \ref{est} for  $\tilde q =r'q$ we have
$$ \int_H \int_H (h(T,x,y))^{\tilde q} \nu (dx) \nu (dy) \le \int_H \int_H h_{\tilde q}(T,x,y) \nu (dx) \nu (dy)$$
$$\le \int_H \int_H k_{\tilde q} \exp \{k_{\tilde q} (|x| + \int_0^T |B_3(s)y|ds )\}\nu (dx) \nu (dy)$$
$$\le k_{\tilde q} \int_H \exp \{k_{\tilde q} |x|\} \nu (dx) \int_H \exp \{ k_{\tilde q} \int_0^T |B_3(s)y|ds )\} \nu (dy)$$
$$= k_{\tilde q} \mathbb{E} e^{k_{\tilde q} |\tilde Z|}\cdot \mathbb{E} \exp \{k_{\tilde q} \int_0^T |B_3(s)\tilde Z|ds \} $$
where $\tilde Z$ is an arbitrary random variable with probability distribution $\nu $. By (\ref{y5}), (\ref{SFP}) and the Fernique inequality we
conclude that (\ref{817}) holds true. The proof of (\ref{817}) for the case when $\frac{q}{p'}<1$ is even simpler and is omitted.
\begin{remark}\label{iny}
There is a natural question whether the transition density is regular (continuous) "in $y$", that is, whether the
mapping $y \mapsto d(T,x,y)$ is continuous, at least on a certain subspace $\hat H \subset H$) of full measure. In the
Gaussian case the formulas for the density may be used to conclude that if $\overline{S_T^*Q^{-1}_T} \in
\mathcal{L}(\hat H,H)$ then $y\to g(T,x,y)$ is continuous on $\hat H$ for all $T>0$ and $x\in H$ (cf. the
Cameron-Martin formula (\ref{y1a})). A similar well-known formula for $k(T,y)$ (see e.g. \cite{reg}) yields $\hat H\to
H$ continuity of the mapping $y\mapsto k(T,y)$ provided
\begin{equation}\label{y13}
C(T) := \overline{Q_{\infty}^{-1/2}(I-S_0(T)S^*_0(T))^{-1}S_0(T)S_0^*(T)Q_{\infty}^{-1/2}} \in \mathcal{L}(\hat H,H)
\end{equation}
where $S_0(T) = Q_{\infty}^{-1/2}S_TQ_{\infty}^{1/2}$. Following the proof of Theorem \ref{dens} we can easily see that
the remaining factor, the function $h(T,x,y)$ is continuous in $y\in\hat{H}$ if the mapping $y \to \hat Z^{x,y}_t$ is
$\hat H \to H$ a.s. continuous (which by Proposition \ref{cont} (b) happens if $\overline{K_tQ_T^{-1/2}} \in
\mathcal{L}(\hat H,H)$) and
\begin{equation}\label{B3}
B_3 \in \mathcal{L}(\hat H, L^1(0,T;H)).
\end{equation}
We are able to verify these additional conditions in some important cases (supposing that the standing assumptions of this Section (\ref{GIR}), (\ref{h4}) and (\ref{ea1})
are satisfied).
\par\noindent
(a) All three conditions are satisfied if $\mathrm{dim}{H}<\infty$.
\par\noindent
(b) In the commutative case the first two conditions are satisfied with $\hat{H}=H$. However, condition
(\ref{B3}) is not
satisfied with $\hat H=H$ even in simple infinite - dimensional situations and a smaller space $\hat H$ must be considered (cf. Example
\ref{example} below for details).
\par\noindent
\end{remark}
\begin{example}\label{example}
Consider the semilinear stochastic heat equation
\begin{equation}\label{HEQ}
\frac{\partial u}{\partial t}(t,\xi ) = \frac{\partial ^2 u}{\partial \xi ^2} (t,\xi ) + f(u(t,\xi )) + \eta (t,\xi ),
\quad (t, \xi )\in \mathbb{R}_+ \times (0,1),
\end{equation}
with an initial condition and Dirichlet boundary conditions
\begin{equation}\label{HEQ2}
u(0, \xi)= x(\xi ),\quad u(t,0)=u(t,1)=0,\quad t\ge 0,\, \xi \in (0,1)
\end{equation}
where $f: \mathbb{R} \to \mathbb{R}$ is bounded and continuous and $\eta $ denoted formally a space-dependent white
noise. As well known (see e.g. \cite{dz1} for fundamentals on the theory of stochastic evolution equations) the system
(\ref{HEQ}) - (\ref{HEQ2}) may be understood as an equation of the form (\ref{SEM}) in the space $H = L^2(0,1)$ where
$A= \frac{\partial ^2}{\partial \xi ^2}$, $\mathrm{dom} (A) = H^1_0(0,1)\cap H^2(0,1)$, $F: H \to H$, $F(y)(\xi) :=
f(y(\xi ))$, $ y\in H$, $\xi \in (0,1)$, and $\sqrt{Q}$ is a bounded operator on $H=L^2(0,1)$. We assume that the
operator $Q$ is boundedly invertible on $H$, (i.e., the noise is nondegenerate). Then Hypothesis \ref{GIR} is obviously
satisfied and Hypothesis \ref{h4} (ii) is satisfied with $\beta = \frac{1}{2}$ and arbitrary $\alpha \in (0,
\frac{1}{2})$ (cf.\cite{den}, Example 9.2 and references therein). Thus the conclusions of Theorems \ref{dens},
\ref{HSS} and \ref{lp} hold true in the present example. \vskip1cm\noindent
 As far as continuity of the transition
density "in the variable $y$" is concerned (cf. Remark \ref{iny} ), the problem is more difficult and we only can
verify our conditions in the diagonal (commutative) case. Denote by $(e_n)$ and $(\alpha _n)$ the orthonormal basis in
$H$ consisting of eigenvectors of the operator $-A$ and its corresponding eigenvalues (so we have $\alpha _n >0$, $
\alpha _n \sim n^2$), and assume that $Q$ commutes with $A$, that is,
$$ Qe_n = \lambda _n e_n, \quad 0< \inf \lambda _n \le \sup \lambda _n <\infty. $$
Then it is easy to compute eigenvalue expansions of all operators that are needed in Remark \ref{iny}. We have
\begin{equation}\label{y14}
K_tQ_T^{-1}e_n = \frac{1-e^{-2\alpha _n t}}{1-e^{-2\alpha _n T}}e^{-\alpha _n (T-t)} e_n,
\end{equation}
\begin{equation}\label{y15}
Q_T^{-1}S^*_T e_n= 2e^{-\alpha _n T} \frac{\alpha _n}{\lambda _n}(1-e^{-\alpha _n T})^{-1} e_n,
\end{equation}
\begin{equation}\label{y15a}
C(T) e_n= 2e^{-2\alpha _n T} \frac{\alpha _n}{\lambda _n}(1-e^{-2\alpha _n T})^{-1} e_n,
\end{equation}
\begin{equation}\label{y16}
B_3(s) e_n = 2e^{-\alpha _n (T-s)} \frac{\alpha _n}{\sqrt{\lambda _n}}(1-e^{-\alpha _n T})^{-1} e_n.
\end{equation}
Obviously, all operators given in (\ref{y14}) and (\ref{y15a}) are in $\mathcal{L} (H)$, but it is easy to see that
$\| B_3(s) \| \sim \frac{1}{T-s}$, so $B_3$ is not an element of $\mathcal{L} (H,L^1(0,T;H))$
 and we do not obtain the continuity in $y$ in the norm of $H$. However,
taking $\hat H = \mathrm{dom} ((-A)^\delta )$ endowed with the graph norm for
any $\delta >0$ (which coincides with a suitable Sobolev-Slobodetskii space) we may
easily check that the condition (\ref{B3}) is satisfied and we may conclude that the mapping $y\mapsto d(T,x,y)$ is $\hat H\to H$ continuous.
\par\noindent
In the present case it is also easy to write equation (\ref{b2}) for the OU Bridge that splits into a sequence of
independent one-dimensional equations for particular coordinates $\hat z_n^{x,y}(t) := \left\langle \hat
Z^{x,y}_t,e_n\right\rangle$. We obtain
$$d\hat z_n^{x,y}(t) = [-\alpha _n\hat z_n^{x,y}(t)-2\alpha_ne^{-\alpha_n(T-t)}(1-e^{-2\alpha_n (T-t)})^{-1}
(e^{-\alpha_n(T-t)}\hat z_n^{x,y}(t) -y_n)]dt + \sqrt{\lambda_n}d \zeta_n(t) $$
for $t\in(0,T)$ with the initial condition
$$\hat z_n^{x,y}(0) = x_n,$$
where $x_n = \left\langle x, e_n \right\rangle$,  $y_n = \left\langle y, e_n \right\rangle$ and $\zeta _n(t) =
\left\langle \zeta _t, e_n \right\rangle $. Here we do not have to assume that the eigenvalues $\alpha_n$ are all
negative, only $\alpha_n \neq 0$. If $\alpha_n =0$ for some $n$ the corresponding equation takes the form
$$d\hat z_n^{x,y}(t) = \frac{y_n-\hat z_n^{x,y}(t)}{T-t}dt + \sqrt{\lambda_n}d \zeta_n(t),\quad t\in(0,T),$$
which is a well-known equation for a one-dimensional Brownian Bridge.
\end{example}


\begin{thebibliography}{99}
\bibitem{ball}
Ball J. M.: Strongly continuous semigroups, weak solutions, and the variation of constants formula \emph{Proc. Amer. Math. Soc.}  63  (1977), 370-373
\bibitem{ania-mild}
Chojnowska-Michalik A.: Stochastic differential equations in Hilbert spaces, in: Probability theory (Papers, VIIth Semester, Stefan Banach Internat. Math. Center, Warsaw, 1976),  pp. 53-74, Banach Center Publ., 5, PWN, Warsaw, 1979
\bibitem{ania}
Chojnowska-Michalik A.: Transition semigroups for stochastic semilinear equations on Hilbert spaces   \emph{Dissertationes Math.} 396  (2001)
\bibitem{fock}
Chojnowska-Michalik A. and Goldys B.: Nonsymmetric Ornstein-Uhlenbeck semigroup as second quantized operator.  \emph{J. Math. Kyoto Univ.} 36  (1996), 481-498
\bibitem{reg}
Chojnowska-Michalik A. and Goldys B.: On regularity properties of nonsymmetric Ornstein-Uhlenbeck semigroup in $L^p$ spaces,  \emph{Stochastics and
Stochastics Rep.} 59  (1996), 183-209
\bibitem{symm}
Chojnowska-Michalik A. and Goldys B.: Symmetric
Ornstein-Uhlenbeck Semigroups and their Generators,
{\em Probab. Theory and Related Fields\/} 124 (2002), 459-486
\bibitem{dz1}
Da Prato G. and Zabczyk J.: Stochastic Equations in Infinite Dimensions, Cambridge University Press 1992
\bibitem{furman}
Fuhrman M.: Regularity properties of transition probabilities in infinite dimensions \emph{Stochastics and Stochastics Rep.}  69  (2000), 31-65
\bibitem{anal}
Goldys B.: On analyticity of Ornstein-Uhlenbeck semigroups, \emph{Atti Accad. Naz. Lincei Cl. Sci. Fis. Mat. Natur. Rend. Lincei} (9) Mat. Appl.  10  (1999), 131-140
\bibitem{acta}
Goldys B. and van Neerven J.M.A.M.: Transition semigroups of Banach space valued Ornstein-Uhlenbeck
processes, {\it Acta Appl. Math.} 76 (2003), 283-330
\bibitem{den}
Goldys B. and Maslowski B.: Lower estimates of transition densities and bounds on exponential ergodicity for stochastic
PDE's, \emph{Ann. Probab.} 34 (2006), 1451-1496
\bibitem{lyons}
Lyons T. J and Zheng W.A.: On Conditional Diffusion Processes,
{\em Proc. Royal Soc. Edinburgh\/} 115A (1990), 243-255
\bibitem{ma}
Ma Zhi Ming and R\"ockner M.: Introduction to the theory of (nonsymmetric) Dirichlet forms, Springer-Verlag, 1992
\bibitem{mandelbaum}                                                                                    Mandelbaum A.: Linear estimators and measurable linear transformations on a Hilbert space, {\it Z. Wahrsch. Verw. Gebiete} 65 (1984), 385-397
\bibitem{masi1}
Maslowski B. and Sim\~ao I.: Asymptotic properties of
stochastic semilinear equations by the method of lower
measures, {\em Colloquium Math.\/} 72 (1997), 147-171
\bibitem{masi2}
Maslowski B. and Sim\~ao I.: Long time behaviour of
non-autonomous SPDE's, {\em Stochastic Processes and }
{\em Applications\/} 95 (2001), 285-309
\bibitem{simao1}
Sim\~ao I.: Pinned Ornstein-Uhlenbeck processes on an infinite-dimensional space, Stochastic Analysis and Applications (Powys, 1995), World Sci. Publishing, River Edge, NJ, 1996.
\bibitem{neerven}
van Neerven J.M.A.M.: Nonsymmetric Ornstein-Uhlenbeck
Semigroups in Banach Spaces, {\em J. Funct. Anal.\/} 155 (1998),
495-535
\bibitem{on}
Ondrej\' at M.: Brownian representations of cylindrical martingales, martingale problem and strong Markov property of weak solutions of SPDEs in Banach spaces, {\em Czechoslovak Math. J.\/} 55 (2005), 1003-1039
\bibitem{yor}
Yor M.: Some Aspects of Brownian Motion, Birkh\"auser 1992
\end{thebibliography}
\end{document}